\documentclass[lettersize,journal]{IEEEtran}
\usepackage{amsmath,amsfonts,amssymb}
\usepackage{algorithmic}
\usepackage{algorithm}
\usepackage{array}
\usepackage[caption=false,font=normalsize,labelfont=sf,textfont=sf]{subfig}
\usepackage{textcomp}
\usepackage{stfloats}
\usepackage{url}
\usepackage{verbatim}
\usepackage{graphicx}
\usepackage{cite}
\usepackage{booktabs}
\usepackage{xcolor}
\newtheorem{proposition}{Proposition}

\begin{document}

\setlength{\abovecaptionskip}{3pt}
\setlength{\textfloatsep}{0.8\baselineskip plus 0.2\baselineskip minus 0.4\baselineskip}
\setlength{\floatsep}{0.6\baselineskip plus 0.2\baselineskip minus 0.2\baselineskip}
\makeatletter
\def\@IEEENORMtitlevspace{1\baselineskip}
\def\@IEEEMINtitlevspace{0.75\baselineskip}
\makeatother

\title{Data Center Power Flexibility Characterization: A DeepSets Reinforcement Learning Approach}

\author{Qiong~Liu,
        Changhong~Zhao,~\IEEEmembership{Senior Member,~IEEE,}
        Meng~Yang,
        Xi~Weng,
        and~Yue~Chen,~\IEEEmembership{Senior Member,~IEEE}% <-this % stops a space
\thanks{This work was supported by the
National Natural Science Foundation of China (NSFC) and the Hong Kong Research Grants Council (RGC)  under the
NSFC/RGC Joint Research Scheme (Grant No. 72561160161 and N\_CUHK422/25), and the Hong Kong RGC under Grant C1049-24G. (Corresponding author: Yue Chen.)}% <-this % stops a space
\thanks{Q. Liu, M. Yang, and Y. Chen are with the Department of Mechanical and Automation Engineering, The Chinese University of Hong Kong, Hong Kong, China (email: liuqiong\_yl@outlook.com; myang@mae.cuhk.edu.hk; yuechen@mae.cuhk.edu.hk).}% <-this % stops a space
\thanks{C. Zhao is with the Department of Information Engineering, The Chinese University of Hong Kong, Hong Kong, China (email: chzhao@ie.cuhk.edu.hk).}
\thanks{X. Weng is with the Guanghua School of Management, Peking University, China (email: wengxi125@gsm.pku.edu.cn).}}

% The paper headers
\markboth{}%
{Liu \MakeLowercase{\textit{et al.}}: Data Center Power Flexibility Characterization: A DeepSets Reinforcement Learning Approach}

\maketitle

\begin{abstract}
Data centers can provide power flexibility through the scheduling of computing workloads, cooling systems, and batteries. Characterizing this flexibility online is useful but challenging because of the complex coupling among these resources.
To address this challenge, this paper proposes a role-conditioned cross-head DeepSets reinforcement learning approach to characterize the power flexibility interval of a data center at each time step. First, we formulate the flexibility characterization problem as a Markov decision process that jointly models the three resources under a carbon budget, with actions specifying the upper and lower endpoints of the power flexibility interval. We then design DeepSets actor and critic that represent delay-tolerant workloads as an unordered set, enabling the proposed approach to accommodate varying numbers of workloads with a fixed number of trainable parameters. Finally, because maximizing the width of the flexibility interval requires moving its endpoints in opposite directions, we develop a role-conditioned cross-head decoder that distinguishes their respective roles while keeping them coordinated.
Case studies on the Google Cluster Data show that the proposed approach provides more than twice the daily flexibility of the standard soft actor-critic method with nearly zero carbon overshoot while using 87.7\% fewer trainable parameters.

%Data centers can provide power flexibility through the scheduling of computing workloads, cooling devices, and batteries. Characterizing this flexibility online is useful but challenging because of the complex coupling among these resources.
%To address this challenge, this paper proposes a role-conditioned cross-head DeepSets reinforcement learning approach, which reports a power flexibility interval to the system operator at each time step. We formulate the flexibility characterization as a Markov decision process that models the three resources jointly under a carbon budget and whose action specifies both endpoints of the interval. Since the number of delay-tolerant workloads varies over time, we design a DeepSets actor and critic that treat these workloads as an unordered set and accommodate any number of workloads with a fixed number of trainable parameters. Since maximizing the flexibility interval pushes its two endpoints in opposite directions, we develop a role-conditioned cross-head decoder that distinguishes their roles while keeping them coordinated.
%Case studies on the Google Cluster Data show that the proposed approach provides more than twice the daily flexibility of the standard soft actor-critic method with nearly zero carbon overshoot while using 87.7\% fewer trainable parameters, and this reduction grows as the number of workloads increases.
% When the number of workloads increases 25-fold, it attains about 19 times the return of the soft actor-critic with 98.7\% fewer parameters.
\end{abstract}

\begin{IEEEkeywords}
Carbon-aware operation, data center, deep reinforcement learning, power flexibility.
\end{IEEEkeywords}

\section{Introduction}
\IEEEPARstart{D}{ata} centers are among the fastest-growing electricity loads, with their electricity consumption rising from about 415~TWh in 2024 to a projected 945~TWh by 2030, mainly driven by artificial intelligence~\cite{ieaEnergyAI2025}. Rather than acting as passive consumers, data centers can provide flexibility to the power system through their computing workloads, cooling systems, and co-located battery energy storage systems (BESSs)~\cite{chenInternetDataCenters2020, caoDatadrivenFlexibilityAssessment2022, hanEvaluatingDispatchableCapacity2026}.
The power flexibility of these three kinds of resources cannot be assessed independently, since workload scheduling changes the IT power and thus the heat dissipated by the cooling systems, and the BESSs adjust the total grid power through charging and discharging. Meanwhile, the rising electricity demand of data centers increases carbon emissions~\cite{radovanovicCarbonAwareComputing2023, hallCarbonAwareComputing2025}, making emissions a significant consideration in their operation.

The energy management of data centers has been extensively studied, whereas few studies have quantified their power flexibility. Existing energy management studies mainly focus on workload scheduling~\cite{radovanovicCarbonAwareComputing2023, hallCarbonAwareComputing2025, wangMultiObjectiveLow2026}, power capping~\cite{sunLearningEnabledAdaptivePower2025}, demand response and ancillary services~\cite{chenInternetDataCenters2020, cuiCoOptimizationCoincident2026}, and the coordination of workload scheduling with cooling system control or BESSs~\cite{panUncertaintyawareReinforcementLearning2026, liuMixtureofexpertsBasedMulticritic2026a}. 
Regarding flexibility quantification, the temporal flexibility of periodic batch workloads was assessed from production traces in~\cite{caoDatadrivenFlexibilityAssessment2022}, and the flexibility provided by IT workload shifting and chiller setpoint adjustment in response to carbon intensity was evaluated in~\cite{misaghianAssessmentCarbonAware2023}. The spatial-temporal flexibility of cross-regional data center clusters was quantified in~\cite{hanEvaluatingDispatchableCapacity2026}, where backup energy storage and workload uncertainty were considered. These studies on quantifying data center power flexibility consider only one or two of the three resources (workloads, cooling systems, and energy storage), and they are performed offline using historical traces or workload scenarios.

Since future workload arrivals, carbon intensities, and dispatch signals are unknown in advance, an online power flexibility characterization method is more practical. For easier integration into system operation, online methods typically report a power interval at each time step, within which the system operator can dispatch any point without knowing the detailed internal operation of the data center~\cite{yanRealTimeFeedback2024, wangCarbonAwareQuantification2025}. Such intervals were characterized by model predictive control for distribution systems~\cite{chenAggregatePowerFlexibility2020} and buildings~\cite{gasserPredictiveEnergyManagement2021}, and by Lyapunov optimization for electric vehicles (EVs)~\cite{yanRealTimeFeedback2024}, with carbon-aware virtual queues incorporated in~\cite{wangCarbonAwareQuantification2025} to keep the carbon emissions within a budget. However, these model-based online methods may not be suitable for quantifying data center power flexibility due to the large number of delay-tolerant workloads and their complex coupling with cooling systems and BESSs.

Learning-based methods provide another approach to online flexibility characterization, in which a policy is learned from operational data without requiring an explicit system model. To date, such methods have been developed mainly for EVs, where real-time aggregate flexibility feedback was learned in~\cite{liLearningBasedPredictive2021} and a soft actor-critic (SAC) agent was trained in~\cite{zhangLearningBasedRealTime2024} to provide the base power and aggregate flexibility of a charging station. For data centers, deep reinforcement learning (DRL) has been applied mainly to the energy management reviewed above, including workload scheduling~\cite{wangMultiObjectiveLow2026}, power capping~\cite{sunLearningEnabledAdaptivePower2025}, joint workload, cooling, and BESS control~\cite{naugSustainDCBenchmarkingSustainable2024a, liuMixtureofexpertsBasedMulticritic2026a}, and market bidding~\cite{panUncertaintyawareReinforcementLearning2026}. 
However, the learning-based flexibility characterization of a data center has not been well explored, and directly applying the methods above is challenging. This is because the number of active delay-tolerant workloads varies over time and can be large, resulting in high-dimensional state and action spaces whose dimensions change dynamically in DRL. Standard actor and critic networks have fixed input and output dimensions, making it difficult to accommodate such large, varying dimensional spaces. A graph neural network was adopted in~\cite{orfanoudakisScalableReinforcementLearning2025} to represent a varying number of EVs on a graph of chargers and vehicles, but the delay-tolerant workloads in a data center are not connected by a graph topology and form an unordered set.

Besides the large and varying number of active workloads, the flexibility interval itself poses a further challenge. The DRL policy is trained to maximize the interval width under the operational constraints, so the two endpoints of the interval play opposite roles, as the upper endpoint is driven upward and the lower endpoint downward. The policy must distinguish the two roles while keeping the endpoints coordinated, since the reward depends on their difference. The flexibility feedback learned in~\cite{liLearningBasedPredictive2021} was a probability distribution over the aggregate power levels rather than an interval. In~\cite{zhangLearningBasedRealTime2024}, both bounds of the flexibility region were derived from the aggregate output of a single actor. Neither of these methods has addressed the challenge posed by the two endpoints.
% However, the challenge that the two endpoints pose was not considered in these methods. 

To fill the research gaps above, this paper develops a role-conditioned cross-head (RCCH) DeepSets reinforcement learning approach for online carbon-aware power flexibility characterization of a data center. Our main contributions are:
\begin{itemize}
    \item We formulate the online carbon-aware power flexibility characterization for a data center as a Markov decision process with a two-branch action that specifies the two endpoints of the power flexibility interval, and design a constraint-aware action map that converts each action into executable controls. Any requested power within the interval is disaggregated across the individual workloads, the cooling system, and the BESS by interpolating the two endpoint controls with a single dispatch ratio, so that the system operator does not need to know the detailed internal operation of the data center.
    \item We design a DeepSets actor and critic that treat the active workloads as an unordered set, so that the cooling and BESS actions and the critic value are permutation-invariant with respect to the workloads, whereas the workload actions are permutation-equivariant. Compared to a standard SAC whose input and output dimensions must be fixed in advance, this design handles a varying number of active workloads, reduces the number of trainable parameters, and improves the learning performance.
    \item We propose an RCCH decoder tailored to the opposing optimization directions of the two endpoints of the power flexibility interval, one driven upward and the other downward. It represents the actions of the two branches by role-specific tokens to distinguish their roles, and mixes the tokens across the two branches through a self-attention mechanism to keep them coordinated, which further improves the learning performance.
\end{itemize}

The rest of this paper is organized as follows. Section~\ref{sec:2} presents the model for quantifying the power flexibility of data centers. Section~\ref{sec:3} develops the RCCH DeepSets reinforcement learning approach. Section~\ref{sec:results} presents the case studies, and Section~\ref{sec:conclusion} concludes the paper.

\section{Carbon-Aware Power Flexibility Formulation}
\label{sec:2}
\label{sec:2.flex}

\begin{figure}[t]
    \centering
    \includegraphics[width=\columnwidth]{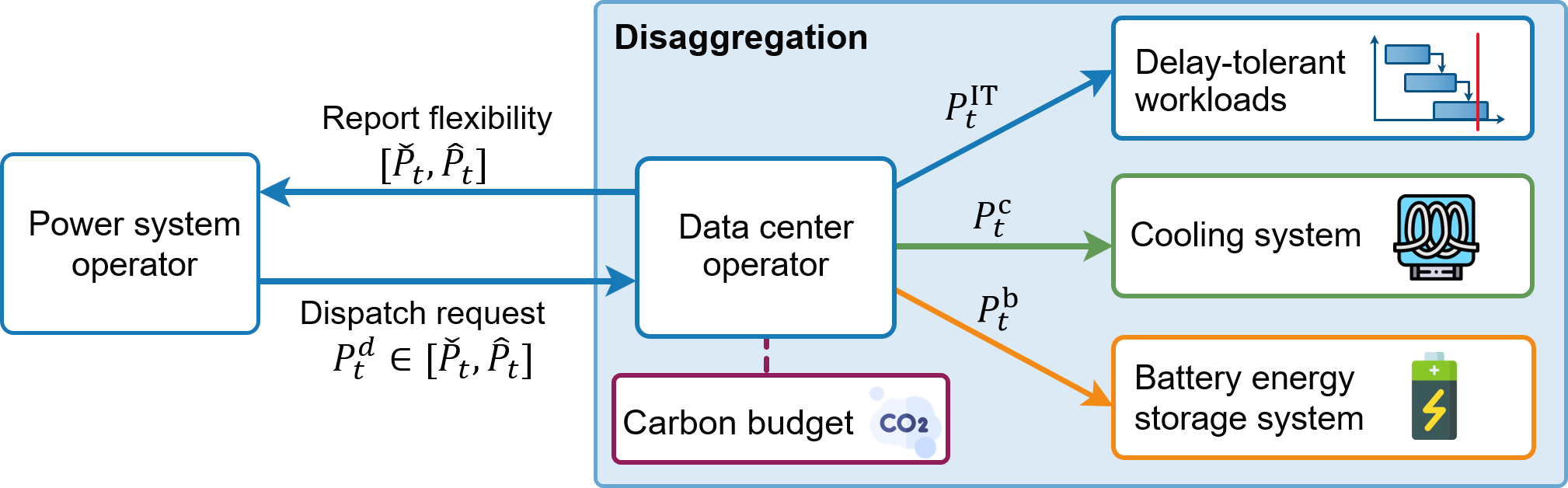}
    \caption{System diagram of data center power flexibility characterization.}
    \label{fig:system_diagram}
\end{figure}

The data center derives power flexibility from three kinds of coupled resources: delay-tolerant workloads can be scheduled at any time step before their deadlines, the thermal inertia of the facility allows cooling demand to be shifted over time, and the BESS shifts grid power through charging and discharging.
Over a horizon of $T$ time steps of length $\Delta t$, indexed by $t\in\mathcal{T}=\{1,\dots,T\}$, the data center reports a power flexibility interval $[\check P_t,\hat P_t]$ to the system operator at each time step $t$, and the system operator then requests a power $P_t^{\mathrm d}\in[\check P_t,\hat P_t]$, which the data center disaggregates to the workloads, the cooling system, and the BESS%
% , as shown in Fig.~\ref{fig:system_diagram}
. Since the requests may increase the grid power or shift it toward time steps of high carbon intensity, the carbon emissions over the horizon are further required to remain within a carbon budget. The overall procedure is shown in Fig. \ref{fig:system_diagram}.

In this section, we first establish an offline model to characterize data center power flexibility, assuming complete knowledge of the future workloads and external conditions.
%(i.e., ambient temperature and grid carbon intensity). 
To make it more practical, we then introduce an online feedback-based procedure to solve it.

\subsection{Offline Model}
\label{sec:2.offline}

The offline model jointly optimizes an upper and a lower endpoint trajectory, whose grid power $\hat P_t$ and $\check P_t$ bound the interval at each time step, to maximize the total interval width:
\begin{equation}\label{eq:offline_band_obj}
 \mathbf{P}^{\mathrm{off}}:\quad
 \max_{\hat{\mathbf{u}}_{1:T},\,\check{\mathbf{u}}_{1:T}}\;\sum\nolimits_{t\in\mathcal{T}}\big(\hat{P}_t-\check{P}_t\big),
\end{equation}
where $\hat{\mathbf{u}}_t=\big(\{\hat\mu_{j,t}\}_{j\in\mathcal{J}},\,\hat P_t^{\mathrm{c}},\,\hat P_t^{\mathrm{ch}},\,\hat P_t^{\mathrm{dis}}\big)$ collects the upper-endpoint controls, namely, the service rates $\hat\mu_{j,t}$, the cooling power $\hat P_t^{\mathrm{c}}$, and the BESS charging and discharging powers $\hat P_t^{\mathrm{ch}}$ and $\hat P_t^{\mathrm{dis}}$, whose difference $\hat P_t^{\mathrm b}=\hat P_t^{\mathrm{ch}}-\hat P_t^{\mathrm{dis}}$ is the net BESS power. The vector $\check{\mathbf{u}}_t$ collects the lower-endpoint controls.

The computing workloads are either latency-sensitive or delay-tolerant. The latency-sensitive workload $L_t$, expressed as a fraction of the total computing capacity, is served upon arrival, whereas each delay-tolerant workload $j\in\mathcal{J}$ is characterized by its arrival time step $a_j$, last executable time step $d_j$, required workload $W_j$, and maximum service rate $\mu_j^{\max}$~\cite{liModelingDemandResponse2015}. The service rate $\hat\mu_{j,t}$ is the fraction of the capacity allocated to workload $j$ at time step $t$. The constraints of the upper trajectory are given by
\begingroup\allowdisplaybreaks
\begin{subequations}\label{eq:off_con}
\begin{align}
 & \hat w_{j,t+1}=\hat w_{j,t}-\hat\mu_{j,t}\Delta t,~\forall j,\forall t\in[a_j,d_j], \label{eq:off_work}\\
 &\hat w_{j,a_j}=W_j,~\hat w_{j,d_j+1}=0,~\forall j, \label{eq:offline_band_complete}\\
 & 0\leq\hat\mu_{j,t}\leq\mu_j^{\max},~\forall j,\forall t\in[a_j,d_j], \label{eq:off_rate}\\
 & \hat\mu_{j,t}=0,~\forall j,\forall t\notin[a_j,d_j], \label{eq:off_zero}\\
& \sum\nolimits_{j\in\mathcal{J}}\hat\mu_{j,t}\leq\bar\mu_t,~\forall t, \label{eq:off_cap}\\
& \hat P_t^{\mathrm{IT}}=P_{\mathrm{static}}+P_{\mathrm{dyn}}\big(L_t+\sum\nolimits_{j\in\mathcal{J}}\hat\mu_{j,t}\big),~\forall t, \label{eq:off_itpower}\\
 & \hat\theta_{t+1}=\hat\theta_t+\tfrac{\Delta t}{C^{\mathrm{th}}}\big(\hat P_t^{\mathrm{IT}}-\kappa_t\hat P_t^{\mathrm c} -U(\hat\theta_t-\theta_t^{\mathrm{amb}})\big),~\forall t, \label{eq:off_thermal}\\
 & 0\leq\hat P_t^{\mathrm c}\leq P^{\mathrm{c,max}},~\forall t, \label{eq:off_cbound}\\
 & \theta^{\min}\leq\hat\theta_{t+1}\leq\theta^{\max},~\forall t, \label{eq:off_tbound}\\
 & \hat e_{t+1}=\hat e_t+\big(\eta^{\mathrm{ch}}\hat P_t^{\mathrm{ch}}-\hat P_t^{\mathrm{dis}}/\eta^{\mathrm{dis}}\big)\Delta t,~\forall t, \label{eq:off_battery}\\
 & e^{\min}\leq\hat e_{t+1}\leq e^{\max},~\forall t, \label{eq:off_soc}\\
 & 0\leq\hat P_t^{\mathrm{ch}}\leq P^{\mathrm{ch,max}},~0\leq\hat P_t^{\mathrm{dis}}\leq P^{\mathrm{dis,max}},~\forall t, \label{eq:off_rating}\\
 & \hat P_t^{\mathrm{ch}}\hat P_t^{\mathrm{dis}}=0,~\forall t, \label{eq:off_comp}\\
 & \hat P_t =\hat P_t^{\mathrm{IT}}+\hat P_t^{\mathrm c}+\hat P_t^{\mathrm{ch}}-\hat P_t^{\mathrm{dis}},~\forall t, \label{eq:off_grid}
\end{align}
\end{subequations}
\endgroup
where $\hat w_{j,t}$, $\hat\theta_t$, and $\hat e_t$ are the remaining workload, hall temperature, and BESS energy on the upper trajectory, respectively. Constraints~\eqref{eq:off_work} and~\eqref{eq:offline_band_complete} describe the remaining workload and require each workload to be completed by its last executable time step. Constraints~\eqref{eq:off_rate} and~\eqref{eq:off_zero} limit the service rate and restrict service to $[a_j,d_j]$. Constraint~\eqref{eq:off_cap} limits the total service rate by the available computing capacity $\bar\mu_t=1-L_t$. Constraint~\eqref{eq:off_itpower} gives the upper-endpoint IT power $\hat P_t^{\mathrm{IT}}$, where $P_{\mathrm{static}}$ and $P_{\mathrm{dyn}}$ are the static IT power and the dynamic IT power at full utilization. Constraint~\eqref{eq:off_thermal} gives the hall temperature dynamics, where $C^{\mathrm{th}}$, $U$, $\theta_t^{\mathrm{amb}}$, and $\kappa_t$ are the thermal capacitance, heat transfer coefficient, ambient temperature, and cooling coefficient of performance, respectively. Constraints~\eqref{eq:off_cbound} and~\eqref{eq:off_tbound} impose the cooling power limit $P^{\mathrm{c,max}}$ and the admissible hall temperature band $[\theta^{\min},\theta^{\max}]$. Constraints~\eqref{eq:off_battery}--\eqref{eq:off_comp} describe the BESS, where $\eta^{\mathrm{ch}}$ and $\eta^{\mathrm{dis}}$ are the charging and discharging efficiencies, and $[e^{\min},e^{\max}]$, $P^{\mathrm{ch,max}}$, and $P^{\mathrm{dis,max}}$ are the energy limits and power ratings. Constraint~\eqref{eq:off_grid} gives the power purchased from the main grid; let $\hat P_t^{\mathrm b}=\hat P_t^{\mathrm{ch}}-\hat P_t^{\mathrm{dis}}$ denote the net BESS power.

Moreover, the carbon emissions must remain within the budget:
\begin{equation}
 \sum\nolimits_{t\in\mathcal{T}}g_t\hat P_t\Delta t\leq C^{\mathrm{bud}}, \label{eq:offline_band_carbon}
\end{equation}
where $g_t\geq0$ and $C^{\mathrm{bud}}$ denote the grid carbon intensity and the horizon carbon budget, respectively.

Similarly, the lower-endpoint controls $\check{\mathbf{u}}_t$ and the corresponding states $\check w_{j,t}$, $\check\theta_t$, and $\check e_t$ satisfy the same constraints, and both trajectories start from the same initial hall temperature $\theta_1$ and BESS energy $e_1$.

The two trajectories are further coupled by the following constraints, which guarantee that every requested power sequence within the intervals is achievable, as stated in Proposition~\ref{prop-1}:
\begin{align}
 \hat\mu_{j,t}&\geq\check\mu_{j,t},~\forall j,\forall t, \label{eq:off_share}\\
 \kappa_t\hat P_t^{\mathrm c}-\hat P_t^{\mathrm{IT}}&\geq\kappa_t\check P_t^{\mathrm c}-\check P_t^{\mathrm{IT}},~\hat P_t^{\mathrm b}\geq\check P_t^{\mathrm b},~\forall t. \label{eq:off_order}
\end{align}
These constraints also imply $\hat P_t\geq\check P_t$, so each power flexibility interval is well defined.

\begin{proposition}\label{prop-1}
    For any requested power sequence within the intervals obtained by the offline model, i.e., $P_t^{\mathrm d}\in[\check P_t,\hat P_t]$ for all $t$, there exists an achievable carbon-aware control strategy.
\end{proposition}

The proof of Proposition~\ref{prop-1} is given in Appendix~\ref{app:prop1}. Despite this guarantee, the offline model requires accurate forecasts of the future workloads and external conditions. Furthermore, since the requests are unknown when the intervals are determined, the intervals must remain achievable for every request sequence, which may underestimate the power flexibility (e.g., that for delay-tolerant workloads). This motivates us to develop an online feedback-based flexibility characterization method that can better unlock data center power flexibility by taking into account the latest information.%In addition, the model fails to capture the flexibility of the delay-tolerant workloads, since constraints~\eqref{eq:offline_band_complete} and~\eqref{eq:off_share} require both trajectories to complete each workload $W_j$ with $\hat\mu_{j,t}\geq\check\mu_{j,t}$ at every time step, which forces the service rates on the two trajectories to coincide.

\subsection{Online Feedback-Based Flexibility Characterization}
\label{sec:2.dispatch_aware}

The online flexibility characterization determines the intervals sequentially based on the current state without forecasts, and uses each request of the system operator as feedback for the subsequent intervals.
At each time step $t$: 1) the data center observes the current state and external conditions, determines the upper- and lower-endpoint controls $\hat{\mathbf{u}}_t$ and $\check{\mathbf{u}}_t$, and reports the resulting power flexibility interval $[\check P_t,\hat P_t]$ to the system operator; 2) the system operator requests a power $P_t^{\mathrm d}\in[\check P_t,\hat P_t]$; 3) the data center disaggregates the request across its computing, cooling, and BESS resources; and 4) the data center updates the remaining workloads, hall temperature, BESS energy, and remaining carbon budget $C_t^{\mathrm{rem}}$. The remaining carbon budget is initialized as $C_1^{\mathrm{rem}}=C^{\mathrm{bud}}$ and updated as $C_{t+1}^{\mathrm{rem}}=C_t^{\mathrm{rem}}-g_tP_t^{\mathrm d}\Delta t$.

In step~3), the data center first computes the dispatch ratio $\beta_t=(P_t^{\mathrm d}-\check P_t)/(\hat P_t-\check P_t)$, with $\beta_t=0$ if $\hat P_t=\check P_t$, and then interpolates the endpoint controls with $\beta_t$ to obtain the dispatched quantities of each resource, marked by the superscript $\mathrm d$:
\begin{align}
 \mu_{j,t}^{\mathrm d}
 &= (1-\beta_t)\check\mu_{j,t}
    +\beta_t\hat\mu_{j,t},~\forall j\in\mathcal{J}_t, \nonumber\\
 P_t^{\mathrm{c,d}}
 &= (1-\beta_t)\check P_t^{\mathrm{c}}
    +\beta_t \hat P_t^{\mathrm{c}}, \nonumber\\
 P_t^{\mathrm{b,d}}
 &= (1-\beta_t)\check P_t^{\mathrm{b}}
    +\beta_t \hat P_t^{\mathrm{b}},\label{eq:online_control_interpolation}\\
 P_t^{\mathrm{ch,d}} & =[P_t^{\mathrm{b,d}}]^+,\qquad
 P_t^{\mathrm{dis,d}}=[P_t^{\mathrm{b,d}}]^-,
 \label{eq:online_battery_projection}
\end{align}
where $\mathcal{J}_t\subseteq\mathcal{J}$ is the set of active delay-tolerant workloads at time step $t$, $[x]^+=\max\{x,0\}$, and $[x]^-=\max\{-x,0\}$. The dispatched IT power $P_t^{\mathrm{IT,d}}$ and the realized states then follow from~\eqref{eq:off_itpower}, \eqref{eq:off_work}, \eqref{eq:off_thermal}, and~\eqref{eq:off_battery} with the upper-endpoint variables replaced by the actual dispatched decisions.

\section{RCCH DeepSets Reinforcement Learning}
\label{sec:3}

This section develops the proposed RCCH DeepSets reinforcement learning (RL) approach to achieve the online feedback-based flexibility characterization. We first formulate the characterization as a Markov decision process (MDP) with a two-branch action and design the constraint-aware action map, and then present the DeepSets actor-critic, the RCCH decoder, and the training algorithm.

\subsection{Markov Decision Process Formulation}
\label{sec:3.mdp}

We formulate the online data center power flexibility characterization as an MDP $\mathcal{M}=\langle \mathcal{S},\mathcal{A},\mathcal{P},R,\gamma\rangle$, where $\mathcal{S}$ is the state space, $\mathcal{A}$ is the action space, $\mathcal{P}$ describes the system transition, $R$ is the reward, and $\gamma\in(0,1)$ is the discount factor. The number of active workloads varies over time as workloads arrive and complete, so the dimensions of the state and action vary over time.

\subsubsection{State}
The state $s_t=\big(\{x_{j,t}^{\mathrm{work}}\}_{j\in\mathcal{J}_t},\,x_t^{\mathrm{fac}}\big)$ consists of a workload part and a facility part. Each active workload is described by
\begin{equation}
 x_{j,t}^{\mathrm{work}}=\Big[w_{j,t},\;d_j-t,\;\mu_j^{\max},\;u_{j,t}\Big]^{\top},
\label{eq:method_workload_features}
\end{equation}
where the urgency is $u_{j,t}=\min\{w_{j,t}/(\mu_j^{\max}(d_j-t)\Delta t),\bar{u}\}$ for $t<d_j$ and $u_{j,t}=\bar{u}$ otherwise, with $\bar{u}$ a constant cap. The facility part is
\begin{equation}
\begin{aligned}
 x_t^{\mathrm{fac}}=\big[
 &g_t,\; \cos(2\pi t/T),\; \sin(2\pi t/T),\; \bar{\mu}_t,\\
 &|\mathcal{J}_t|,\; \sum\nolimits_{j\in\mathcal{J}_t}w_{j,t},\; \theta_t^{\mathrm{amb}},\; \theta_t,\\
 &P_{t-1}^{\mathrm{IT,d}},\; P_{t-1}^{\mathrm{c,d}},\; e_t,\; C_t^{\mathrm{rem}}
 \big]^{\top}.
\end{aligned}
\label{eq:method_facility_features}
\end{equation}

\subsubsection{Two-Branch Action}
The action $a_t=\big(a_t^{(1)},a_t^{(2)}\big)\in\mathcal{A}(s_t)$ consists of two branches, where branch $r\in\{1,2\}$ is
\begin{equation}
 a_t^{(r)}=\Big(\{a_{j,t}^{\mu,(r)}\}_{j\in\mathcal{J}_t},
 a_t^{\mathrm{c},(r)},a_t^{\mathrm{b},(r)}\Big)
 \in[-1,1]^{|\mathcal{J}_t|+2},
\label{eq:candidate_action_method}
\end{equation}
so that $\mathcal{A}(s_t)$ has dimension $2|\mathcal{J}_t|+4$. The superscripts $\mu$, $\mathrm{c}$, and $\mathrm{b}$ denote the workload, cooling, and BESS actions.

\subsubsection{Transition}
The next state $s_{t+1}\sim\mathcal{P}(\cdot\,|\,s_t,a_t)$ follows the data center dynamics in~\eqref{eq:off_work}, \eqref{eq:off_thermal}, and~\eqref{eq:off_battery} under the dispatched controls, and is stochastic since $\beta_t$ and the newly arriving workloads are unknown when $a_t$ is selected.

\subsubsection{Reward}
The reward is the sum of a flexibility reward, deadline and thermal penalties, and a carbon term,
\begin{equation}
 r_t=r_t^{\mathrm{p}}+r_t^{\mathrm{ddl}}+r_t^{\theta}+r_t^{\mathrm{cb}}.
\label{eq:method_reward}
\end{equation}
The flexibility reward is
\begin{equation}
 r_t^{\mathrm{p}}=\hat{P}_t-\check{P}_t.
\label{eq:method_reward_flex}
\end{equation}
Since the action map in Section~\ref{sec:3.param} cannot enforce the workload deadlines and the hall temperature band in every state, we introduce the deadline and thermal penalties:
\begin{align}
 r_t^{\mathrm{ddl}}={}&-\lambda^{\mathrm{ddl}}\sum_{j\in\mathcal{J}_t}\big[(t+1)-d_j\big]^+ w_{j,t+1}\Delta t,
\label{eq:method_ddl}\\
 r_t^{\theta}={}&-\lambda^{\theta}\big([\theta_{t+1}-\theta^{\max}]^++[\theta^{\min}-\theta_{t+1}]^+\big).
\label{eq:method_reward_thermal}
\end{align}
The carbon term is
\begin{equation}
\begin{aligned}
 r_t^{\mathrm{cb}}={}&-\lambda^{\mathrm{cb,T}}\,\mathbf{1}_{\{t=T\}}\big[g_TP_T^{\mathrm d}\Delta t-C_T^{\mathrm{rem}}\big]^+\\
 &+\lambda^{\mathrm{cb}}\big(C^{\mathrm{bud}}/(T\Delta t)-g_tP_t^{\mathrm d}\big),
\end{aligned}
\label{eq:method_reward_carbon}
\end{equation}
in which $\mathbf{1}_{\{\cdot\}}$ is the indicator function. The first term penalizes the carbon overshoot at the end of the horizon. Since the terminal penalty is sparse, the second term adds a dense guide that encourages the realized emission rate $g_tP_t^{\mathrm d}$ of each step to stay below the average allowable rate $C^{\mathrm{bud}}/(T\Delta t)$. The deadline, thermal, and terminal carbon penalties are summed over the realized trajectory and the two endpoints. The successor state of each endpoint is obtained by applying its controls to the realized state for one step.

\subsection{Constraint-Aware Action Map}
\label{sec:3.param}

For each endpoint, the constraint-aware action map converts the action into executable controls. The same deterministic maps act on both endpoints, so we omit the endpoint accents.

\subsubsection{Action Saturation}
To widen the interval, the endpoint controls often need to reach the bounds of their admissible ranges, which correspond to the action values $\pm1$. A tanh-squashed Gaussian actor, however, reaches $\pm1$ only asymptotically, and the squashing derivative $1-\tanh^2$ vanishes there. To address this, we rescale each action component by $1/\rho$ and clip it to $[-1,1]$, where $\rho$ is the action saturation threshold, chosen slightly below 1 in the proposed approach:
\begin{equation}
 \tilde{a}=\operatorname{clip}\!\left(a/\rho,\,-1,\,1\right),~\rho\in(0,1].
\label{eq:method_saturation}
\end{equation}

\subsubsection{Workload Map with Deadline Floors}
For each active workload, the map first builds the service range $[\mu_{j,t}^{\mathrm{req}},\mu_{j,t}^{\mathrm{cap}}]$:
\begin{equation}
\begin{aligned}
 \mu_{j,t}^{\mathrm{req}}&=\min\big\{\mu_j^{\max},\big[w_{j,t}/\Delta t-[d_j-t]^+\mu_j^{\max}\big]^+\big\},\\
 \mu_{j,t}^{\mathrm{cap}}&=\min\big\{\mu_j^{\max},w_{j,t}/\Delta t\big\},
\end{aligned}
\label{eq:method_required_service}
\end{equation}
where the deadline floor $\mu_{j,t}^{\mathrm{req}}$ is the least service at time step $t$ that still allows the remaining workload to be completed by its deadline, and the service cap $\mu_{j,t}^{\mathrm{cap}}$ is the largest service the workload can receive.

The workload action is then mapped into this range to obtain the proposed rate $\mu_{j,t}^{\mathrm{prop}}$,
\begin{equation}
 \mu_{j,t}^{\mathrm{prop}}
 =\mu_{j,t}^{\mathrm{req}}+\frac{\tilde{a}_{j,t}^{\mu}+1}{2}
 \big(\mu_{j,t}^{\mathrm{cap}}-\mu_{j,t}^{\mathrm{req}}\big).
\label{eq:method_service_map}
\end{equation}
Each range is built for one workload in isolation, so these proposed rates may jointly exceed the available capacity $\bar{\mu}_t$. The workload map therefore rescales them so that the total rate does not exceed $\bar{\mu}_t$. With the totals $S_t^{\mathrm{req}}=\sum\nolimits_{j\in\mathcal{J}_t}\mu_{j,t}^{\mathrm{req}}$ and $S_t^{\mathrm{prop}}=\sum\nolimits_{j\in\mathcal{J}_t}\mu_{j,t}^{\mathrm{prop}}$, the rate of each $j\in\mathcal{J}_t$ is
\begin{equation}
 \mu_{j,t}=
 \renewcommand{\arraystretch}{1.25}
 \left\{\begin{array}{@{}l@{\;}l@{}}
 \mu_{j,t}^{\mathrm{prop}}, & S_t^{\mathrm{prop}}\le \bar{\mu}_t,\\
 \mu_{j,t}^{\mathrm{req}}+\zeta_t(\mu_{j,t}^{\mathrm{prop}}-\mu_{j,t}^{\mathrm{req}}), & S_t^{\mathrm{req}}\le \bar{\mu}_t<S_t^{\mathrm{prop}},\\
 \mu_{j,t}^{\mathrm{EDF}}, & S_t^{\mathrm{req}}>\bar{\mu}_t,
 \end{array}\right.
\label{eq:method_capacity_projection}
\end{equation}
where the common scaling factor of the second case is
\begin{equation}
 \zeta_t=\frac{\bar{\mu}_t-S_t^{\mathrm{req}}}{S_t^{\mathrm{prop}}-S_t^{\mathrm{req}}}.
\label{eq:method_capacity_scale}
\end{equation}
In the second case, each workload keeps its deadline floor, and only the part of its proposed rate above the floor is scaled by $\zeta_t$. In the third case, the deadline floors exceed $\bar{\mu}_t$ and cannot all be met; the map then serves them in earliest-deadline-first (EDF) order until $\bar{\mu}_t$ is exhausted. Let $\sigma_t(n)$ denote the $n$th workload in $\mathcal{J}_t$ sorted by nondecreasing deadline. For $n=1,\ldots,|\mathcal{J}_t|$, the EDF rate in \eqref{eq:method_capacity_projection} is
\begin{equation}
 \mu_{\sigma_t(n),t}^{\mathrm{EDF}}=\min\big\{\mu_{\sigma_t(n),t}^{\mathrm{req}},\big[\bar{\mu}_t-\sum\nolimits_{\ell=1}^{n-1}\mu_{\sigma_t(\ell),t}^{\mathrm{EDF}}\big]^+\big\},
\label{eq:method_edf_fallback}
\end{equation}
where $\sum\nolimits_{\ell=1}^{n-1}\mu_{\sigma_t(\ell),t}^{\mathrm{EDF}}$ is the rate already allocated to the $n-1$ workloads with earlier deadlines.

\subsubsection{Cooling Map}
The cooling map keeps the next hall temperature within the band $[\theta^{\min},\theta^{\max}]$ whenever this is achievable with a cooling power in $[0,P^{\mathrm{c,max}}]$. By~\eqref{eq:off_thermal}, the cooling power that brings the next hall temperature to a target $\vartheta$, clipped to $[0,P^{\mathrm{c,max}}]$, is
\begin{equation}
\begin{aligned}
 \bar{P}_t^{\mathrm{c}}(\vartheta)=\min\big\{&P^{\mathrm{c,max}},\big[\big(P_t^{\mathrm{IT}}-U(\theta_t-\theta_t^{\mathrm{amb}})\\
 &-C^{\mathrm{th}}(\vartheta-\theta_t)/\Delta t\big)/\kappa_t\big]^+\big\}.
\end{aligned}
\label{eq:cooling_limit_map}
\end{equation}
A lower target temperature requires more cooling power, so the map sets the cooling power range to $P_t^{\mathrm{c,req}}=\bar{P}_t^{\mathrm{c}}(\theta^{\max})$ and $P_t^{\mathrm{c,cap}}=\bar{P}_t^{\mathrm{c}}(\theta^{\min})$, and maps the cooling action affinely into this range,
\begin{equation}
 P_t^{\mathrm{c}}=P_{t}^{\mathrm{c,req}}
 +\frac{\tilde{a}_t^{\mathrm{c}}+1}{2}
 \big(P_{t}^{\mathrm{c,cap}}-P_{t}^{\mathrm{c,req}}\big).
\label{eq:cooling_affine_map}
\end{equation}

\subsubsection{BESS Map}
Given the current BESS energy $e_t$, the charging and discharging powers are limited to
\begin{equation}
\begin{aligned}
 \bar{P}_t^{\mathrm{ch}}&=\min\big\{P^{\mathrm{ch,max}},\big[(e^{\max}-e_t)/(\eta^{\mathrm{ch}}\Delta t)\big]^+\big\},\\
 \bar{P}_t^{\mathrm{dis}}&=\min\big\{P^{\mathrm{dis,max}},\big[\eta^{\mathrm{dis}}(e_t-e^{\min})/\Delta t\big]^+\big\},
\end{aligned}
\label{eq:method_battery_limits}
\end{equation}
so the net BESS power lies in the range $[-\bar{P}_t^{\mathrm{dis}},\bar{P}_t^{\mathrm{ch}}]$. The BESS action is then mapped onto this range by
\begin{equation}
 P_t^{\mathrm{b}}=
 \begin{cases}
 \tilde{a}_t^{\mathrm{b}}\bar{P}_t^{\mathrm{ch}}, & \tilde{a}_t^{\mathrm{b}}\geq 0,\\
 \tilde{a}_t^{\mathrm{b}}\bar{P}_t^{\mathrm{dis}}, & \tilde{a}_t^{\mathrm{b}}<0.
 \end{cases}
\label{eq:method_battery_map}
\end{equation}

\subsubsection{Endpoint Ordering Map}
This map ensures that each workload rate, the cooling power, and the BESS power of the upper endpoint are no smaller than those of the lower endpoint. The saturated actions of the two branches are first ordered componentwise, $\check{a}_t=\min\{\tilde a_t^{(1)},\tilde a_t^{(2)}\}$ and $\hat{a}_t=\max\{\tilde a_t^{(1)},\tilde a_t^{(2)}\}$. The workload map converts $\check{a}_t$ into the lower rates $\check\mu_{j,t}$. Since the rescaling in~\eqref{eq:method_capacity_projection} may reduce an individual rate obtained from $\hat{a}_t$ below its lower counterpart, the upper rates are instead defined relative to the lower rates:
\begin{equation}
 \hat\mu_{j,t}=\check\mu_{j,t}+\xi_t\big[\mu_{j,t}'-\check\mu_{j,t}\big]^+,~\forall j\in\mathcal{J}_t,
\label{eq:method_order_rates}
\end{equation}
where $\mu_{j,t}'$ are the rates given by the workload map for $\hat{a}_t$ and $\xi_t\in[0,1]$ is the largest factor that keeps the total upper rate within $\bar\mu_t$. The IT powers $\check P_t^{\mathrm{IT}}$ and $\hat P_t^{\mathrm{IT}}$ are then computed from $\check\mu_{j,t}$ and $\hat\mu_{j,t}$ using~\eqref{eq:off_itpower}, and the cooling and BESS maps convert $\check{a}_t$ and $\hat{a}_t$ into the corresponding controls of the two endpoints.

\subsection{DeepSets Actor-Critic}
\label{sec:3.net}

The standard actor-critic approach, built on a multilayer perceptron (MLP) whose input concatenates the facility features and the features of all workloads into one vector, faces two challenges. First, the number of active workloads $|\mathcal{J}_t|$ varies over time, whereas the input widths of the actor and critic and the output width of the actor are fixed by a padding bound $J^{\max}$, set to an estimated maximum number of active workloads. If $|\mathcal{J}_t|<J^{\max}$, the input is zero-padded to $J^{\max}$ entries, whereas if $|\mathcal{J}_t|>J^{\max}$, the excess workloads cannot be represented. Second, since a data center can run many delay-tolerant workloads at once, $J^{\max}$ must be large, so the input layers of the actor and critic and the output layer of the actor contain a large number of parameters, which may make them more difficult to train.

In addition, the active workloads have no natural order, so listing them in a different order describes the same state. Accordingly, the network outputs should satisfy two properties: 1) permutation invariance of the cooling and BESS actions and the critic value, which describe the data center as a whole; and 2) permutation equivariance of the workload actions, which are defined per workload. For example, if workloads~1 and~2 swap positions in the input, the cooling and BESS actions remain unchanged, while the workload actions of workloads~1 and~2 swap positions in the output. The MLP assigns separate weights to each input position, so it has neither property by construction.

\begin{figure}[t]
    \centering
    \includegraphics[width=0.8\columnwidth]{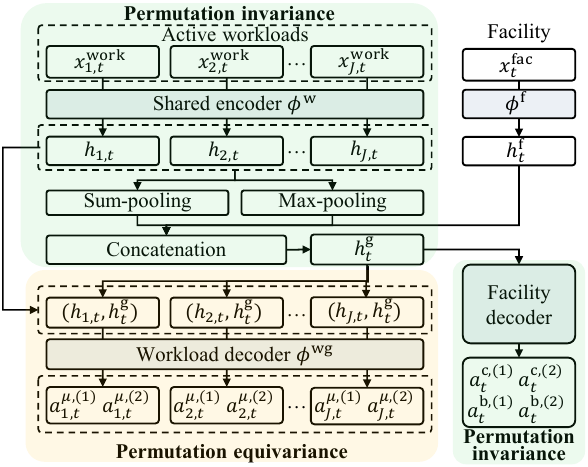}
    \caption{Architecture of the DeepSets actor with $J=|\mathcal{J}_t|$ active workloads.}
    \label{fig:deepsets_encoder}
\end{figure}

To address these challenges and satisfy the two properties, we build on the DeepSets construction~\cite{zaheerDeepSets2017} and design a DeepSets actor and critic. The structure of the actor is shown in Fig.~\ref{fig:deepsets_encoder}. A shared encoder first maps every active workload to an embedding,
\begin{equation}
 h_{j,t}=\phi^{\mathrm{w}}(x_{j,t}^{\mathrm{work}}).
\label{eq:method_workload_encoder}
\end{equation}
The embeddings are then summarized by scaled sum-pooling and max-pooling,
\begin{equation}
 c_t^{\mathrm{sum}}
 =\frac{1}{N^{\mathrm{ref}}}\sum_{j\in\mathcal{J}_t}h_{j,t},
 \qquad
 c_t^{\max}
 =\max_{j\in\mathcal{J}_t}h_{j,t},
\label{eq:method_deepsets_pooling}
\end{equation}
in which the pooling scale $N^{\mathrm{ref}}$ is a fixed reference number of workloads.

The pooled summary is concatenated with the encoded facility features into a global embedding:
\begin{equation}
 \begin{aligned}
 h_t^{\mathrm{f}}&=\phi^{\mathrm{f}}(x_t^{\mathrm{fac}}),\\
 h_t^{\mathrm{g}}&=\big(h_t^{\mathrm{f}},\,c_t^{\mathrm{sum}},\,c_t^{\max}\big).
 \end{aligned}
\label{eq:method_facility_embedding}
\end{equation}
This global embedding is combined with each workload embedding by a workload decoder $\phi^{\mathrm{wg}}$, which produces the two actions of that workload,
\begin{equation}
 \big(a_{j,t}^{\mu,(1)},a_{j,t}^{\mu,(2)}\big)
 =\phi^{\mathrm{wg}}\big(h_{j,t},\,h_t^{\mathrm{g}}\big),
\label{eq:method_workload_decoder}
\end{equation}
and a facility decoder converts $h_t^{\mathrm{g}}$ into the two pairs of cooling and BESS actions $\{(a_t^{\mathrm{c},(r)},a_t^{\mathrm{b},(r)})\}_{r=1}^{2}$. Since both poolings are independent of the workload order, the cooling and BESS actions are permutation-invariant, and sharing $\phi^{\mathrm{wg}}$ across workloads makes the workload actions permutation-equivariant.

The critic appends the two workload actions to $x_{j,t}^{\mathrm{work}}$ and the four cooling and BESS actions to $x_t^{\mathrm{fac}}$, and estimates $Q(s_t,a_t)$ by another DeepSets encoder of the same form followed by an MLP head. Since the pooling is independent of the workload order, the Q-value is permutation-invariant.

\subsection{Role-Conditioned Cross-Head Decoder for the Endpoints}
\label{sec:3.rcch}

\begin{figure}[t]
    \centering    \includegraphics[width=0.51\columnwidth]{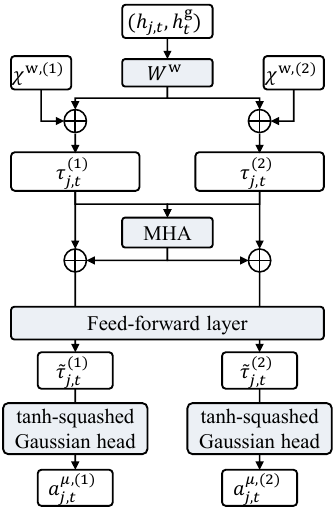}
    \caption{Workload branch of the role-conditioned cross-head decoder.}
    \label{fig:rcch_decoder}
\end{figure}

Maximizing the interval width $\hat{P}_t-\check{P}_t$ drives the two action branches in opposite directions, and the reward depends on their difference. A decoder shared by both branches, such as $\phi^{\mathrm{wg}}$ in~\eqref{eq:method_workload_decoder}, does not represent their roles explicitly. We therefore propose an RCCH decoder that distinguishes the two action branches while keeping them coordinated, as shown in Fig.~\ref{fig:rcch_decoder}.

Each action is represented by an action token that combines the encoded features with a learned embedding of its role. In the workload branch, the two action tokens of workload $j$ are
\begin{equation}
 \tau_{j,t}^{(r)}=W^{\mathrm{w}}\big(h_{j,t},\,h_t^{\mathrm{g}}\big)+\chi^{\mathrm{w},(r)},
 \qquad r\in\{1,2\},
\label{eq:method_rcch_job_tokens}
\end{equation}
where $W^{\mathrm{w}}$ projects the concatenated workload and global embeddings into a token space of dimension $d_{\mathrm{tok}}$, and $\chi^{\mathrm{w},(r)}\in\mathbb{R}^{d_{\mathrm{tok}}}$ is a learned role embedding. The two tokens are mixed by a multi-head self-attention (MHA) layer~\cite{vaswaniAttentionAllYou2017}:
\begin{equation}
 \tau_{j,t}^{\prime(r)}=\tau_{j,t}^{(r)}+\big[\mathrm{MHA}\big(\mathrm{LN}(\tau_{j,t}^{(1)}),\mathrm{LN}(\tau_{j,t}^{(2)})\big)\big]_r,
\label{eq:method_rcch_mixer}
\end{equation}
in which $\mathrm{LN}$ is layer normalization~\cite{baLayerNormalization2016}, which normalizes each token to zero mean and unit variance over its entries and then applies a learned elementwise scale and shift, and $[\cdot]_r$ denotes its $r$th output. Through the attention layer, each token receives information from the token of the other branch. A feed-forward network $\mathrm{FFN}$, which consists of two linear layers with a hidden dimension of $2d_{\mathrm{tok}}$ and a sigmoid-weighted linear unit (SiLU) activation~\cite{elfwingSigmoidweightedLinearUnits2018}, is applied to each $\tau_{j,t}^{\prime(r)}$ separately:
\begin{equation}
 \tilde{\tau}_{j,t}^{(r)}=\mathrm{LN}\Big(\tau_{j,t}^{\prime(r)}+\mathrm{FFN}\big(\mathrm{LN}(\tau_{j,t}^{\prime(r)})\big)\Big),
 \qquad r\in\{1,2\}.
\label{eq:method_rcch_ffn}
\end{equation}
For each branch $r$, a linear output layer maps $\tilde{\tau}_{j,t}^{(r)}$ to the mean and log standard deviation of a Gaussian, and the Gaussian sample is squashed by tanh to give $a_{j,t}^{\mu,(r)}$.

The facility branch uses a similar architecture, with the global embedding $h_t^{\mathrm{g}}$ as input and the cooling and BESS actions as output. Each token carries two role embeddings, one for the action branch $r$ and one for the action type $m$:
\begin{equation}
 \tau_{t}^{(r,m)}=W^{\mathrm{f}}h_t^{\mathrm{g}}+\chi^{\mathrm{f},(r)}+\chi^{\mathrm{f},(m)},
 \qquad m\in\{\mathrm{c},\mathrm{b}\},
\label{eq:method_rcch_fac_tokens}
\end{equation}
where $W^{\mathrm{f}}$ projects $h_t^{\mathrm{g}}$ into the same token space and $\chi^{\mathrm{f},(r)},\chi^{\mathrm{f},(m)}\in\mathbb{R}^{d_{\mathrm{tok}}}$ are learned. The four tokens $\{\tau_t^{(r,m)}\}$ are then processed by MHA and feed-forward layers.

\subsection{Training Algorithm}
\label{sec:3.train}

The proposed training algorithm builds on SAC~\cite{haarnojaSoftActorCriticOffPolicy2018}. The stochastic DeepSets actor with the RCCH decoder samples $a_t\sim\pi_{\omega}(\cdot|s_t)$, where $\omega$ collects all actor parameters, and the replay buffer $\mathcal{D}$ stores the transitions $(s_t,a_t,r_t,s_{t+1})$. The DeepSets critic is instantiated as two Q-networks $Q_{\psi_1}$ and $Q_{\psi_2}$, which are trained by minimizing
\begin{equation}
 \mathcal{L}_{Q_i}=
 \mathbb{E}_{\mathcal{D}}\left[
 \left(Q_{\psi_i}(s_t,a_t)-y_t\right)^2
 \right],\qquad i=1,2,
\label{eq:method_critic_loss}
\end{equation}
with the soft Bellman target
\begin{align}
 y_t={}&r_t+\gamma\bigg(\min_{i=1,2}Q_{\bar{\psi}_i}(s_{t+1},a_{t+1}) \notag\\
 &-\alpha\log\pi_{\omega}(a_{t+1}|s_{t+1})\bigg),
\label{eq:method_critic_target}
\end{align}
where $a_{t+1}\sim\pi_{\omega}(\cdot|s_{t+1})$, $\alpha$ is the entropy temperature, and $Q_{\bar{\psi}_i}$ are the target Q-networks. The actor loss, with $s_t\sim\mathcal{D}$ and $a_t\sim\pi_{\omega}(\cdot|s_t)$, is
\begin{equation}
 \mathcal{L}_{\pi}
 =\mathbb{E}_{\mathcal{D},\pi_{\omega}}\left[\alpha\log\pi_{\omega}(a_t|s_t)
 -\min_{i=1,2}Q_{\psi_i}(s_t,a_t)\right].
\label{eq:method_actor_loss}
\end{equation}
The temperature $\alpha$ is adjusted automatically to keep the policy entropy close to a target value $\mathcal{H}_t^{\mathrm{targ}}$ by minimizing, under the same sampling,
\begin{equation}
 \mathcal{L}_{\alpha}
 =\mathbb{E}_{\mathcal{D},\pi_{\omega}}\left[
 -\alpha\left(
 \log\pi_{\omega}(a_t|s_t)
 +\mathcal{H}_{t}^{\mathrm{targ}}
 \right)
 \right].
\label{eq:method_temperature_loss}
\end{equation}

Because $|\mathcal{J}_t|$ changes over time, each sample has its own target entropy, set to the negative dimension of $\mathcal{A}(s_t)$,
\begin{equation}
 \mathcal{H}_{t}^{\mathrm{targ}}=-\big(2|\mathcal{J}_t|+4\big).
\label{eq:method_dynamic_entropy}
\end{equation}
Algorithm~\ref{alg:method_training} summarizes the proposed training procedure.

\begin{algorithm}[t]
    \caption{Training algorithm of the proposed approach.}
\label{alg:method_training}
\begin{algorithmic}[1]
\REQUIRE Actor $\pi_{\omega}$, critic Q-networks $Q_{\psi_1},Q_{\psi_2}$, target Q-networks $Q_{\bar{\psi}_1},Q_{\bar{\psi}_2}$, replay buffer $\mathcal{D}$, warm-up steps $K_0$, total training steps $K$
\STATE Initialize actor, critic, target, and temperature parameters; set $k\leftarrow 0$
\WHILE{$k<K$}
    \STATE Reset the environment and observe $s_1$
    \FOR{$t=1,2,\ldots,T$}
        \IF{$k<K_0$}
            \STATE Sample $a_t$ uniformly from $\mathcal{A}(s_t)$
        \ELSE
            \STATE Sample the two-branch action $a_t\sim\pi_{\omega}(\cdot|s_t)$
        \ENDIF
        \STATE Apply the constraint-aware action map to obtain the endpoint controls
        \STATE Report $[\check{P}_t,\hat{P}_t]$; sample $\beta_t$ and execute the dispatched controls given by~\eqref{eq:online_control_interpolation}--\eqref{eq:online_battery_projection}
        \STATE Observe $r_t$ and $s_{t+1}$; store $(s_t,a_t,r_t,s_{t+1})$ in $\mathcal{D}$
        \IF{$k\geq K_0$}
            \STATE Sample a mini-batch from $\mathcal{D}$
            \STATE Update the two Q-networks by minimizing~\eqref{eq:method_critic_loss}
            \STATE Update the actor by minimizing~\eqref{eq:method_actor_loss}
            \STATE Update the entropy temperature by minimizing~\eqref{eq:method_temperature_loss} and softly update the target Q-networks
        \ENDIF
        \STATE $k\leftarrow k+1$
    \ENDFOR
\ENDWHILE
\end{algorithmic}
\end{algorithm}

\section{Case Studies}
\label{sec:results}

\subsection{Simulation Environment}

We conduct numerical experiments to evaluate the proposed approach in a data center adapted from the SustainDC environment~\cite{naugSustainDCBenchmarkingSustainable2024a}. Each daily episode contains $T=96$ steps of 15~min. The ambient temperature of San Jose and the average grid carbon intensity of California are taken from the SustainDC dataset and linearly interpolated from hourly values. The computing demand is derived from the Google Cluster Data~\cite{reissGoogleClusterUsage2011}, whose scheduling classes indicate latency sensitivity; classes 2--3 and 0--1 correspond to the latency-sensitive workload $L_t$ and the aggregate delay-tolerant workload, respectively. Because the trace provides no individual deadlines, delay-tolerant workloads are generated synthetically, with arrivals following a Poisson process calibrated to the class 0--1 workload. The required workload $W_j$, in capacity-hours, is lognormal with mean 0.025 and coefficient of variation 2.0. The maximum service rate $\mu_j^{\max}$ is $W_j$ divided by a minimum service time drawn uniformly from $[0.5,4]$~h, and the deadline $d_j$ is set to the arrival time plus six times this service time, rounded up to a whole step. The IT, cooling, and BESS parameters are listed in Table~\ref{tab:env_parameters}, and the cooling coefficient of performance is set to $\kappa_t=\min\{6.5,\max\{3.0,5.0-0.1(\theta_t^{\mathrm{amb}}-25)\}\}$.

Every tenth day of the year, 37 days in total, is used for testing, and the remaining 328 days are used for training. The dispatch ratio $\beta_t$ is drawn uniformly from $[0,1]$ at each step, i.e., the requested power is sampled uniformly from the reported interval.
The daily carbon budget is set to 95\% of the emissions of the greedy baseline on the same day.
All experiments are implemented in Python~3.10.18 using PyTorch~2.7.1 with CUDA~11.8.

\begin{table}[t]
    \caption{Simulation environment parameters}
    \label{tab:env_parameters}
    \centering
    \footnotesize
    \begin{tabular}{@{}>{\raggedright\arraybackslash}p{\dimexpr0.67\columnwidth-\tabcolsep}>{\raggedright\arraybackslash}p{\dimexpr0.33\columnwidth-\tabcolsep}@{}}
        \toprule
        Parameter & Setting \\
        \midrule
        IT static/dynamic power $P_{\mathrm{static}}$, $P_{\mathrm{dyn}}$ (MW) & 1.10, 2.40 \\
        \midrule
        Thermal capacitance $C^{\mathrm{th}}$ (MWh/$^\circ$C) & 0.10 \\
        Heat transfer coefficient $U$ (MW/$^\circ$C) & 0.02 \\
        Cooling power limit $P^{\mathrm{c,max}}$ (MW) & 1.10 \\
        Hall temperature band $[\theta^{\min},\theta^{\max}]$ ($^\circ$C) & $[18,27]$ \\
        Initial hall temperature ($^\circ$C) & 22 \\
        \midrule
        BESS energy limits $[e^{\min},e^{\max}]$ (MWh) & $[0.2,2.0]$ \\
        BESS charge/discharge power ratings $P^{\mathrm{ch,max}}$, $P^{\mathrm{dis,max}}$ (MW) & 0.5, 0.5 \\
        BESS charge/discharge efficiencies $\eta^{\mathrm{ch}}$, $\eta^{\mathrm{dis}}$ & 0.96, 0.96 \\
        Initial BESS energy (MWh) & 1.0 \\
        \bottomrule
    \end{tabular}
\end{table}

\begin{table}[t]
    \caption{Training and network hyperparameters}
    \label{tab:hyperparameters}
    \centering
    \footnotesize
    \begin{tabular}{@{}>{\raggedright\arraybackslash}p{\dimexpr0.67\columnwidth-\tabcolsep}>{\raggedright\arraybackslash}p{\dimexpr0.33\columnwidth-\tabcolsep}@{}}
        \toprule
        Parameter & Setting \\
        \midrule
        Discount factor $\gamma$ & 0.99 \\
        Replay buffer capacity & $2\times10^{5}$ \\
        Mini-batch size & 256 \\
        Random warm-up steps $K_0$ & 960 \\
        Actor/temperature learning rate & $3\times10^{-4}$ \\
        Critic learning rate & $1\times10^{-3}$ \\
        Target smoothing coefficient & $5\times10^{-3}$ \\
        Initial entropy temperature $\alpha$ & $1\times10^{-3}$ \\
        \midrule
        Deadline penalty coefficient $\lambda^{\mathrm{ddl}}$ & 100 \\
        Thermal penalty coefficient $\lambda^{\theta}$ & 100 \\
        Carbon penalty coefficient $\lambda^{\mathrm{cb,T}}$ & 50 \\
        Carbon guide coefficient $\lambda^{\mathrm{cb}}$ & 2 \\
        \midrule
        Action saturation threshold $\rho$ & 0.9 \\
        Workload encoder widths ($\phi^{\mathrm{w}}$) & 64, 32 \\
        Facility encoder widths ($\phi^{\mathrm{f}}$) & 128, 64 \\
        Critic head hidden dimension & 128 \\
        MLP hidden dimensions & 256, 128 \\
        Padding bound $J^{\max}$ (MLP only) & 256 \\
        Token dimension $d_{\mathrm{tok}}$ & 64 \\
        Attention heads & 4 \\
        Feed-forward hidden dimension & 128 \\
        Pooling scale $N^{\mathrm{ref}}$ & 128 \\
        Urgency cap $\bar{u}$ & 5 \\
        Hidden layer activation & ReLU \\
        \bottomrule
    \end{tabular}
\end{table}

\subsection{Benchmark Methods}

We compare the proposed RCCH DeepSets-SAC (RCCH DS-SAC) approach with three standard RL baselines and two non-learning benchmarks. The RL baselines are SAC~\cite{haarnojaSoftActorCriticOffPolicy2018}, twin delayed deep deterministic policy gradient (TD3)~\cite{fujimotoAddressingFunctionApproximation2018}, and deep deterministic policy gradient (DDPG)~\cite{lillicrapContinuousControlDeep2016}. Their actors and critics are MLPs whose workload inputs are zero-padded to $J^{\max}$ entries. The greedy baseline sets all actions of the two branches to $-1$ and $+1$, respectively, so the action map returns the widest interval of the current step without regard to future steps or the carbon budget. The ideal full-information mixed-integer linear programming (MILP) benchmark solves the online feedback-based flexibility characterization in Section~\ref{sec:2.dispatch_aware} by assuming that the workloads, external conditions, and dispatch ratios $\beta_t$ of the whole test day are known in advance, so that the endpoint controls of all time steps are optimized jointly; its result serves as an ideal reference for other approaches. The returns of both non-learning benchmarks are computed with the same reward as the RL approaches to enable a direct comparison.
Table~\ref{tab:hyperparameters} lists the hyperparameters of the proposed approach; the standard RL baselines use the same hyperparameters whenever applicable, except for the action saturation threshold $\rho=1$. All learned methods use the constraint-aware action map. Every policy is trained for $K=2\times10^{5}$ steps on the training days with three random seeds and evaluated deterministically on all test days every 2,000 steps, with the results averaged over the test days.
For each seed, the final benchmark and ablation results are averaged over the last ten evaluations, and error bars and shaded bands denote the sample standard deviation over the three seeds. The non-learning benchmarks are evaluated once on each test day.

\subsection{Overall Performance Comparison}
\label{sec:results_overall_performance}

\begin{figure}[t]
    \centering
    \includegraphics[width=0.9\columnwidth,trim=0 6.0 0 6.2,clip]{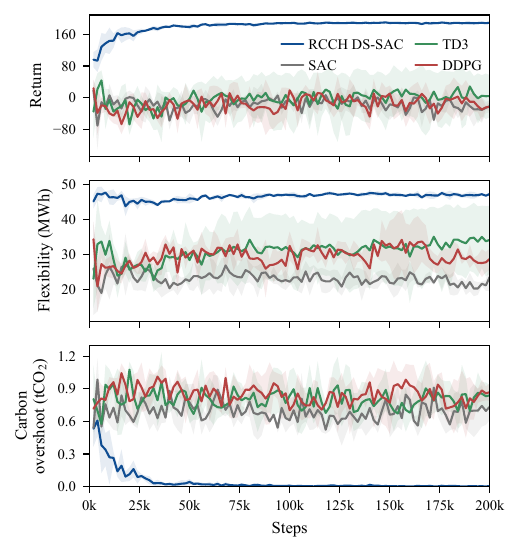}
    \caption{Evaluation results on the test days during training.}
    \label{fig:overall_performance_training}
\end{figure}

As shown in Fig.~\ref{fig:overall_performance_training}, RCCH DS-SAC learns rapidly and attains the highest return and daily flexibility among the learned methods, whereas the returns of the three standard RL baselines fluctuate around zero without a clear improvement. The carbon overshoot of RCCH DS-SAC decreases to a value close to zero within about $3\times10^{4}$ steps, while that of the baselines remains between about 0.5 and 1.0~tCO$_2$. RCCH DS-SAC also exhibits narrower shaded bands across seeds. Deadline and thermal penalties remain zero for all learned methods throughout training and are omitted from the figure.

\begin{figure}[t]
    \centering
    \includegraphics[width=0.9\columnwidth,trim=0 7.6 0 5.7,clip]{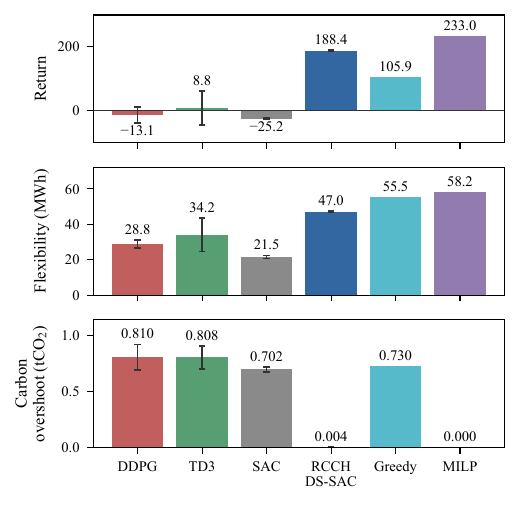}
    \caption{Quantitative comparison of return, daily flexibility, and carbon overshoot on the test days.}
    \label{fig:overall_performance_final}
\end{figure}

The quantitative comparison on the test days is shown in Fig.~\ref{fig:overall_performance_final}. The return of RCCH DS-SAC reaches 188.4, whereas the returns of the three standard RL baselines range from $-25.2$ to 8.8, because their carbon overshoot penalties offset most or all of the flexibility reward. RCCH DS-SAC also provides a daily flexibility of 47.0~MWh, 37.5\% more than TD3, the best-performing RL baseline, and more than twice that of SAC. The greedy baseline provides more flexibility than RCCH DS-SAC but a larger carbon overshoot. The ideal full-information MILP attains the largest flexibility with zero overshoot, but relies on future information that is unavailable in online operation.

As shown in Table~\ref{tab:representation_parameter_comparison}, RCCH DS-SAC contains 87.7\%, 87.1\%, and 80.1\% fewer trainable parameters than SAC, TD3, and DDPG, respectively. Each total comprises one actor and one critic for DDPG and one actor and two critics for the other methods. The MLP baselines are larger because their input layers need a large number of parameters to match the high-dimensional state and action of the data center, whereas the DeepSets encoders of RCCH DS-SAC share their weights across workloads, so that their size does not depend on the number of workloads.

\begin{table}[t]
\caption{Number of trainable network parameters}
\label{tab:representation_parameter_comparison}
\centering
\footnotesize
\setlength{\tabcolsep}{5pt}
\begin{tabular}{lrrrr}
\toprule
 & DDPG & TD3 & SAC & RCCH DS-SAC \\
\midrule
Critic & 430,593 & 430,593 & 430,593 & 29,601 \\
Actor & 364,932 & 364,932 & 431,496 & 99,244 \\
Total & 795,525 & 1,226,118 & 1,292,682 & 158,446 \\
\bottomrule
\end{tabular}
\end{table}

\subsection{Representative Day Flexibility Characterization}
\label{sec:results_representative_day}

\begin{figure}[t]
    \centering
    \includegraphics[width=0.9\columnwidth,trim=0 7.9 0 6.2,clip]{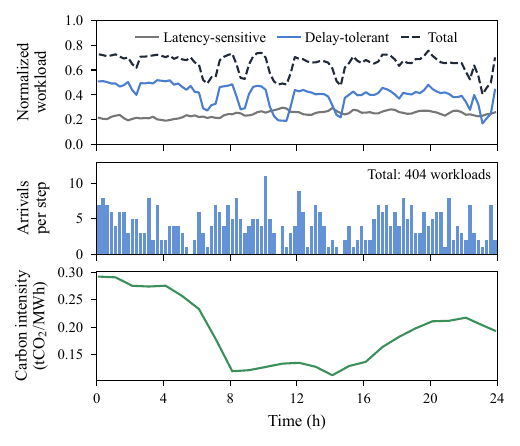}
    \caption{Operating conditions on the representative day.}
    \label{fig:representative_day_conditions}
\end{figure}

\begin{figure}[t]
    \centering
    \includegraphics[width=0.9\columnwidth,trim=0 6.2 0 5.7,clip]{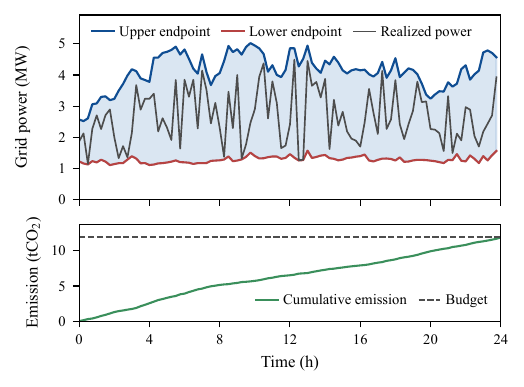}
    \caption{Flexibility interval and realized grid power (top), and cumulative carbon emission with the carbon budget (bottom) on the representative day.}
    \label{fig:representative_day_power_interval}
\end{figure}

We examine the power flexibility characterized by RCCH DS-SAC on a representative test day, the 51st day of the year, whose operating conditions are shown in Fig.~\ref{fig:representative_day_conditions}. As shown in the upper panel of Fig.~\ref{fig:representative_day_power_interval}, the lower endpoint stays between 1.10 and 1.57~MW throughout the day, whereas the upper endpoint ranges from 2.52 to 5.02~MW, so the variation of the interval is mainly due to the upper endpoint. The lower panel shows that the cumulative emission stays within the carbon budget throughout the day.

\begin{figure}[t]
    \centering
    \includegraphics[width=0.9\columnwidth,trim=0 6.4 0 6.4,clip]{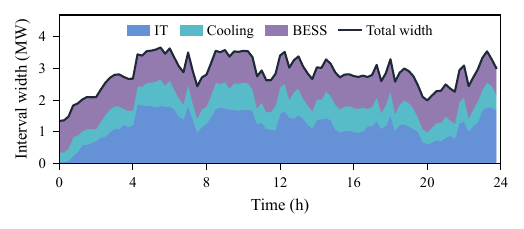}
    \caption{Decomposition of the interval width into the IT, cooling, and BESS components on the representative day.}
    \label{fig:representative_day_endpoint_decomposition}
\end{figure}

To identify the physical sources of this interval, Fig.~\ref{fig:representative_day_endpoint_decomposition} decomposes its width into the IT, cooling, and BESS components. Over the whole day, the IT, cooling, and BESS components contribute 43.0\%, 22.1\%, and 34.9\% of the daily flexibility, respectively. The IT component varies the most, from 0.04 to 1.88~MW, whereas the BESS component stays between 0.82 and 1.00~MW, and the cooling component lies between 0.30 and 0.86~MW.

\begin{figure}[t]
    \centering
    \includegraphics[width=\columnwidth,trim=0 6.4 0 8.8,clip]{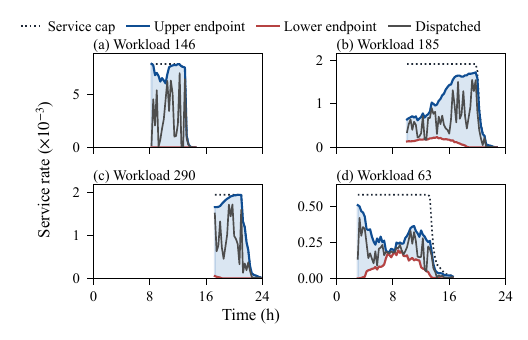}
    \caption{Upper- and lower-endpoint service rates of four delay-tolerant workloads on the representative day.}
    \label{fig:representative_task_flexibility}
\end{figure}

Fig.~\ref{fig:representative_task_flexibility} shows four representative workloads. For workloads 146 and 290, the latter arriving late in the day, the lower endpoint remains close to zero and the upper endpoint closely follows the service cap. For workload 185, the upper endpoint stays well below the service cap after arrival and gradually approaches it near the end of its service. For workload 63, the policy raises the lower endpoint during the morning, bringing it close to the upper endpoint for about two hours.

\subsection{Ablation Study}
\label{sec:results_ablation}

\begin{figure}[!t]
    \centering
    \includegraphics[width=0.9\columnwidth,trim=0 6.9 0 5.7,clip]{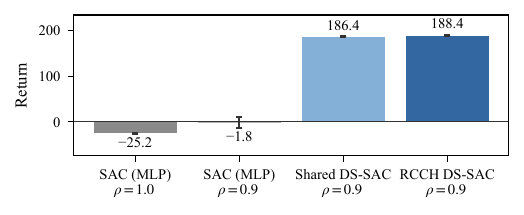}
    \caption{Ablation study of successive techniques and their impact on return.}
    \label{fig:ablation_return_comparison}
\end{figure}

To evaluate the contribution of each proposed technique, we conduct an ablation study, as shown in Fig.~\ref{fig:ablation_return_comparison}. Starting from the standard SAC, we successively introduce the three techniques: applying action saturation with $\rho=0.9$, replacing the padded MLP with the DeepSets workload representation, and replacing the shared decoder with the RCCH decoder. Action saturation increases the mean return, the DeepSets representation yields the largest improvement in return, and the RCCH decoder improves it further.

\subsection{Scalability Study}
\label{sec:results_scalability}

\begin{table}[!t]
\caption{Return and number of trainable parameters at larger numbers of workloads}
\label{tab:cardinality_scalability}
\centering
\footnotesize
\setlength{\tabcolsep}{5pt}
\begin{tabular}{llrrr}
\toprule
 & & & \multicolumn{2}{c}{Parameters} \\
\cmidrule(l){4-5}
Mean $W_j$ & Method & Return & Critic & Actor \\
\midrule
0.01 & SAC & $-18.2$ & 1,610,241 & 1,614,216 \\
 & RCCH DS-SAC & 195.2 & 29,601 & 99,244 \\
\midrule
0.001 & SAC & 8.7 & 15,766,017 & 15,806,856 \\
 & RCCH DS-SAC & 190.9 & 29,601 & 99,244 \\
\bottomrule
\end{tabular}
\end{table}

We further evaluate scalability by increasing the number of delay-tolerant workloads on the one-year dataset. Reducing the mean required workload $W_j$ from its default value of $0.025$ to $0.01$ and $0.001$, with the expected total required workload unchanged, increases the expected number of workloads by $2.5\times$ and $25\times$, respectively. At the two scales, the daily number of delay-tolerant workloads averages 454 and 4,400, respectively, and reaches at most 967 and 9,728. Accordingly, the padding bound $J^{\max}$ of SAC is set to 1,024 and 10,240 at the two scales, and the pooling scale $N^{\mathrm{ref}}$ of RCCH DS-SAC to 320 and 3,200. At each workload scale, RCCH DS-SAC and SAC are trained for $2\times10^{5}$ steps with a single random seed, and the results are averaged over the last ten evaluations.

Table~\ref{tab:cardinality_scalability} reports the resulting return on the test days, together with the trainable parameters in the actor and in each critic, excluding the target networks. RCCH DS-SAC outperforms SAC at both workload scales, with returns above 190, whereas SAC attains $-18.2$ and 8.7 at the $2.5\times$ and $25\times$ scales, respectively. The RCCH actor and critic remain the same size, whereas both SAC networks grow by almost an order of magnitude as the number of workloads increases from the $2.5\times$ to the $25\times$ scale. At the $2.5\times$ and $25\times$ scales, the actor and critics of RCCH DS-SAC together contain 96.7\% and 99.7\% fewer parameters than those of SAC, respectively.

\section{Conclusion}
\label{sec:conclusion}

This paper proposes an RCCH DS-SAC approach for online carbon-aware power flexibility characterization of a data center. We formulate the problem as an MDP whose two-branch action specifies the two endpoints of the flexibility interval, and a constraint-aware action map converts each action into executable controls. Building on this formulation, we develop a DeepSets actor-critic for the varying number of active workloads and an RCCH decoder that distinguishes the roles of the two branches while keeping them coordinated. Case studies on the Google Cluster Data yield the following findings: 1) Among the RL approaches, RCCH DS-SAC attains the highest return and provides more than twice the daily flexibility of SAC with nearly zero carbon overshoot, using 87.7\% fewer trainable parameters. 2) In the ablation study, the DeepSets representation yields the largest improvement in return, and the RCCH decoder improves it further. 3) When the expected number of workloads increases by factors of 2.5 and 25, RCCH DS-SAC outperforms SAC at both scales while keeping its network size constant, whereas the network size of SAC grows by almost an order of magnitude between the two scales. Future work may extend the framework to multiple geographically distributed data centers, where workload placement across sites and locational carbon intensities add a spatial dimension to the flexibility.

\bibliographystyle{IEEEtran}
\bibliography{My_Library}

@article{caoDatadrivenFlexibilityAssessment2022,
  title = {Data-Driven Flexibility Assessment for Internet Data Center towards Periodic Batch Workloads},
  author = {Cao, Yujie and Cheng, Ming and Zhang, Sufang and Mao, Hongju and Wang, Peng and Li, Chao and Feng, Yihui and Ding, Zhaohao},
  year = 2022,
  month = oct,
  journal = {Applied Energy},
  volume = {324},
  note = {{Art.}~no.~119665},
  copyright = {1}
}

@article{chenInternetDataCenters2020,
  title = {Internet Data Centers Participating in Demand Response: {{A}} Comprehensive Review},
  shorttitle = {Internet Data Centers Participating in Demand Response},
  author = {Chen, Min and Gao, Ciwei and Song, Meng and Chen, Songsong and Li, Dezhi and Liu, Qiang},
  year = 2020,
  month = jan,
  journal = {Renewable and Sustainable Energy Reviews},
  volume = {117},
  note = {{Art.}~no.~109466},
  langid = {english}
}

@article{elfwingSigmoidweightedLinearUnits2018,
  title = {Sigmoid-Weighted Linear Units for Neural Network Function Approximation in Reinforcement Learning},
  author = {Elfwing, Stefan and Uchibe, Eiji and Doya, Kenji},
  year = 2018,
  month = nov,
  journal = {Neural Networks},
  volume = {107},
  pages = {3--11},
  doi = {10.1016/j.neunet.2017.12.012},
  langid = {english}
}

@inproceedings{fujimotoAddressingFunctionApproximation2018,
  title = {Addressing Function Approximation Error in Actor-Critic Methods},
  booktitle = {Proceedings of the 35th {{International Conference}} on {{Machine Learning}}},
  author = {Fujimoto, Scott and {van Hoof}, Herke and Meger, David},
  year = 2018,
  month = jul,
  pages = {1587--1596},
  langid = {english}
}

@inproceedings{haarnojaSoftActorCriticOffPolicy2018,
  title = {Soft Actor-Critic: {O}ff-Policy Maximum Entropy Deep Reinforcement Learning with a Stochastic Actor},
  shorttitle = {Soft {{Actor-Critic}}},
  booktitle = {Proceedings of the 35th {{International Conference}} on {{Machine Learning}}},
  author = {Haarnoja, Tuomas and Zhou, Aurick and Abbeel, Pieter and Levine, Sergey},
  year = 2018,
  month = jul,
  pages = {1861--1870},
  langid = {english}
}

@article{hanEvaluatingDispatchableCapacity2026,
  title = {Evaluating the Dispatchable Capacity of Cross-Regional Data Center Clusters Toward Power System Operation},
  author = {Han, Jianpei and Du, Ershun and Du, Bojun and Li, Yaowang and Guo, Jingrong and Zhang, Ning and Kang, Chongqing},
  year = 2026,
  month = mar,
  journal = {IEEE Transactions on Smart Grid},
  volume = {17},
  number = {2},
  pages = {1180--1193},
}

@inproceedings{lillicrapContinuousControlDeep2016,
  title = {Continuous Control with Deep Reinforcement Learning},
  booktitle = {Proceedings of the 4th {{International Conference}} on {{Learning Representations}}},
  author = {Lillicrap, Timothy P. and Hunt, Jonathan J. and Pritzel, Alexander and Heess, Nicolas and Erez, Tom and Tassa, Yuval and Silver, David and Wierstra, Daan},
  year = 2016,
  month = may,
  address = {San Juan, Puerto Rico}
}

@article{liModelingDemandResponse2015,
  title = {Modeling Demand Response Capability by Internet Data Centers Processing Batch Computing Jobs},
  author = {Li, Jie and Bao, Zhen and Li, Zuyi},
  year = 2015,
  month = mar,
  journal = {IEEE Transactions on Smart Grid},
  volume = {6},
  number = {2},
  pages = {737--747}
}

@article{liuMixtureofexpertsBasedMulticritic2026a,
  title = {Mixture-of-Experts Based Multi-Critic Deep Reinforcement Learning for Sustainable Management of Data Center Microgrids},
  author = {Liu, Qiong and Pan, Zhanhua and Guo, Ye and Chen, Yue and Strbac, Goran and Qiu, Dawei},
  year = 2026,
  month = may,
  journal = {Applied Energy},
  volume = {411},
  note = {{Art.}~no.~127561}
}

@inproceedings{naugSustainDCBenchmarkingSustainable2024a,
  title = {{SustainDC}: {B}enchmarking for Sustainable Data Center Control},
  shorttitle = {{{SustainDC}}},
  author = {Naug, Avisek and Guillen, Antonio and Luna, Ricardo and Gundecha, Vineet and Bash, Cullen and Ghorbanpour, Sahand and Mousavi, Sajad and Babu, Ashwin R. and Markovikj, Dejan and Kashyap, Lekhapriya D. and Rengarajan, Desik and Sarkar, Soumyendu},
  year = 2024,
  month = dec,
  booktitle = {Advances in Neural Information Processing Systems},
  volume = {37},
  pages = {100630--100669},
  langid = {english}
}

@techreport{reissGoogleClusterUsage2011,
  author = {Reiss, Charles and Wilkes, John and Hellerstein, Joseph L.},
  title = {{Google} Cluster-Usage Traces: Format + Schema},
  institution = {Google Inc.},
  year = {2011},
  month = nov,
  address = {Mountain View, CA, USA},
  note = {Revised 2014-11-17 for trace version 2.1},
  url = {https://github.com/google/cluster-data}
}

@article{sunLearningEnabledAdaptivePower2025,
  title = {Learning-Enabled Adaptive Power Capping Scheme for Cloud Data Centers},
  author = {Sun, Yimeng and Ding, Zhaohao and Dehghanian, Payman and Teng, Fei},
  year = 2025,
  month = nov,
  journal = {IEEE Transactions on Smart Grid},
  volume = {16},
  number = {6},
  pages = {4755--4767}
}

@inproceedings{zaheerDeepSets2017,
  title = {Deep Sets},
  author = {Zaheer, Manzil and Kottur, Satwik and Ravanbakhsh, Siamak and P{\'o}czos, Barnab{\'a}s and Salakhutdinov, Ruslan and Smola, Alexander J.},
  year = {2017},
  booktitle = {Advances in Neural Information Processing Systems},
  volume = {30}
}

@article{zhangLearningBasedRealTime2024,
  title = {Learning-Based Real-Time Aggregate Flexibility Provision and Scheduling of Electric Vehicles},
  author = {Zhang, Mingyang and Yang, Hongrong and Xu, Yinliang and Sun, Hongbin},
  year = {2024},
  month = nov,
  journal = {IEEE Transactions on Smart Grid},
  volume = {15},
  number = {6},
  pages = {5840--5852}
}

@article{baLayerNormalization2016,
  title = {Layer Normalization},
  author = {Ba, Jimmy Lei and Kiros, Jamie Ryan and Hinton, Geoffrey E.},
  year = {2016},
  journal = {arXiv preprint arXiv:1607.06450}
}

@inproceedings{vaswaniAttentionAllYou2017,
  title = {Attention Is All You Need},
  author = {Vaswani, Ashish and Shazeer, Noam and Parmar, Niki and Uszkoreit, Jakob and Jones, Llion and Gomez, Aidan N. and Kaiser, {\L}ukasz and Polosukhin, Illia},
  year = {2017},
  booktitle = {Advances in Neural Information Processing Systems},
  volume = {30}
}

@article{panUncertaintyawareReinforcementLearning2026,
  title = {Uncertainty-Aware Reinforcement Learning for Joint Optimization of Data Center Operations and Electricity Market Bidding},
  author = {Pan, Zhanhua and Liu, Qiong and Liu, Xiao and Jing, Zhaoxia and Qiu, Dawei},
  year = {2026},
  month = nov,
  journal = {Applied Energy},
  volume = {422},
  note = {{Art.}~no.~128208}
}

@article{radovanovicCarbonAwareComputing2023,
  title={Carbon-Aware Computing for Datacenters},
  volume={38},
  DOI={10.1109/tpwrs.2022.3173250},
  number={2},
  journal={IEEE Transactions on Power Systems},
  author={Radovanovi{\'c}, Ana and Koningstein, Ross and Schneider, Ian and Chen, Bokan and Duarte, Alexandre and Roy, Binz and Xiao, Diyue and Haridasan, Maya and Hung, Patrick and Care, Nick and Talukdar, Saurav and Mullen, Eric and Smith, Kendal and Cottman, MariEllen and Cirne, Walfredo},
  year={2023},
  month=mar,
  pages={1270--1280}
}

@article{hallCarbonAwareComputing2025,
  title={Carbon-Aware Computing for Data Centers with Probabilistic Performance Guarantees},
  DOI={10.1109/tpwrs.2025.3630499},
  journal={IEEE Transactions on Power Systems},
  author={Hall, Sophie and Micheli, Francesco and Belgioioso, Giuseppe and Radovanovi{\'c}, Ana and D{\"o}rfler, Florian},
  year={2025},
  note={early access, doi: 10.1109/TPWRS.2025.3630499}
}

@article{wangMultiObjectiveLow2026,
  title={Multi-Objective Low-Carbon Scheduling Method for Data Centers Based on Ensemble Reinforcement Learning},
  volume={17},
  DOI={10.1109/tsg.2025.3620106},
  number={1},
  journal={IEEE Transactions on Smart Grid},
  author={Wang, Yifan and Sun, Wei and Ren, Pengyu and Harrison, Gareth},
  year={2026},
  month=jan,
  pages={297--308}
}

@article{cuiCoOptimizationCoincident2026,
  title={Co-optimization of Coincident Peak Management and Ancillary Services for Flexible Data Centers},
  DOI={10.1109/tsg.2026.3698492},
  journal={IEEE Transactions on Smart Grid},
  author={Cui, Ruolan and Hu, Chenxi and Yang, Ruizhang and Cui, Xin and Hou, Francis Yunhe},
  year={2026},
  month=nov,
  volume={17},
  number={6},
  pages={5405--5417}
}

@article{misaghianAssessmentCarbonAware2023,
  title={Assessment of Carbon-Aware Flexibility Measures From Data Centres Using Machine Learning},
  volume={59},
  DOI={10.1109/tia.2022.3213637},
  number={1},
  journal={IEEE Transactions on Industry Applications},
  author={Misaghian, M. Saeed and Tardioli, Giovanni and Cabrera, Adalberto Guerra and Salerno, Ilaria and Flynn, Damian and Kerrigan, Ruth},
  year={2023},
  month=jan,
  pages={70--80}
}

@article{chenAggregatePowerFlexibility2020,
  title={Aggregate Power Flexibility in Unbalanced Distribution Systems},
  volume={11},
  DOI={10.1109/tsg.2019.2920991},
  number={1},
  journal={IEEE Transactions on Smart Grid},
  author={Chen, Xin and Dall'Anese, Emiliano and Zhao, Changhong and Li, Na},
  year={2020},
  month=jan,
  pages={258--269}
}

@article{gasserPredictiveEnergyManagement2021,
  title={Predictive energy management of residential buildings while self-reporting flexibility envelope},
  volume={288},
  DOI={10.1016/j.apenergy.2021.116653},
  journal={Applied Energy},
  author={Gasser, Jan and Cai, Hanmin and Karagiannopoulos, Stavros and Heer, Philipp and Hug, Gabriela},
  year={2021},
  month=apr,
  note = {{Art.}~no.~116653}
}

@article{yanRealTimeFeedback2024,
  title={Real-Time Feedback Based Online Aggregate {EV} Power Flexibility Characterization},
  volume={15},
  DOI={10.1109/tste.2023.3324705},
  number={1},
  journal={IEEE Transactions on Sustainable Energy},
  author={Yan, Dongxiang and Huang, Shihan and Chen, Yue},
  year={2024},
  month=jan,
  pages={658--673}
}

@article{wangCarbonAwareQuantification2025,
  title={Carbon-Aware Quantification of Real-Time Aggregate Power Flexibility of Electric Vehicles},
  volume={16},
  DOI={10.1109/tsg.2025.3595577},
  number={6},
  journal={IEEE Transactions on Smart Grid},
  author={Wang, Xiaowei and Chen, Yue},
  year={2025},
  month=nov,
  pages={4906--4920}
}

@article{liLearningBasedPredictive2021,
  title={Learning-Based Predictive Control via Real-Time Aggregate Flexibility},
  volume={12},
  DOI={10.1109/tsg.2021.3094719},
  number={6},
  journal={IEEE Transactions on Smart Grid},
  author={Li, Tongxin and Sun, Bo and Chen, Yue and Ye, Zixin and Low, Steven H. and Wierman, Adam},
  year={2021},
  month=nov,
  pages={4897--4913}
}

@article{orfanoudakisScalableReinforcementLearning2025,
  title={Scalable reinforcement learning for large-scale coordination of electric vehicles using graph neural networks},
  volume={4},
  DOI={10.1038/s44172-025-00457-8},
  number={1},
  note = {{Art.}~no.~118},
  journal={Communications Engineering},
  author={Orfanoudakis, Stavros and Robu, Valentin and Salazar, E. Mauricio and Palensky, Peter and Vergara, Pedro P.},
  year={2025},
  month=jul
}

@techreport{ieaEnergyAI2025,
  title = {Energy and {AI}},
  author = {{International Energy Agency}},
  institution = {International Energy Agency},
  address = {Paris, France},
  year = {2025},
  url = {https://www.iea.org/reports/energy-and-ai}
}

\appendices
\section{Proof of Proposition~1}\label{app:prop1}
\setcounter{equation}{0}
\renewcommand{\theequation}{A.\arabic{equation}}
\begin{IEEEproof}
Suppose that the upper and lower trajectories satisfy~\eqref{eq:off_con}--\eqref{eq:offline_band_carbon} and~\eqref{eq:off_share}--\eqref{eq:off_order}. For any given requested power sequence that satisfies $\check P_t\leq P_t^{\mathrm d}\leq\hat P_t$ for all $t$, we define the dispatch ratio $\beta_t\in[0,1]$ as in Section~\ref{sec:2.dispatch_aware}, so that
\begin{equation}
P_t^{\mathrm d}=(1-\beta_t)\check P_t+\beta_t\hat P_t.
\label{eq:dispatch}
\end{equation}
We claim that the dispatched controls defined by~\eqref{eq:online_control_interpolation}--\eqref{eq:online_battery_projection}, with $\mathcal{J}_t$ replaced by $\mathcal{J}$, form an achievable carbon-aware control strategy for this sequence, i.e., they satisfy~\eqref{eq:off_con}--\eqref{eq:offline_band_carbon} with the upper-endpoint quantities replaced by the dispatched quantities and yield the grid power $P_t^{\mathrm d}$. The realized trajectory starts from the same initial conditions as the endpoint trajectories, namely, $w_{j,a_j}=W_j$ for each workload, the hall temperature $\theta_1$, and the BESS energy $e_1$.

First, we verify the time-decoupled constraints. The dispatched service rates and cooling power are convex combinations of the endpoint values, so they satisfy the rate limits in~\eqref{eq:off_rate} and~\eqref{eq:off_zero}, the capacity limit in~\eqref{eq:off_cap}, and the cooling power limits in~\eqref{eq:off_cbound}. By~\eqref{eq:online_battery_projection}, the dispatched charging and discharging powers are nonnegative and complementary, and $P_t^{\mathrm{ch,d}}-P_t^{\mathrm{dis,d}}=P_t^{\mathrm{b,d}}$. Each endpoint satisfies~\eqref{eq:off_comp}, so $\hat P_t^{\mathrm{ch}}=[\hat P_t^{\mathrm b}]^+$ and $\check P_t^{\mathrm{dis}}=[\check P_t^{\mathrm b}]^-$. Since $P_t^{\mathrm{b,d}}$ lies in $[\check P_t^{\mathrm b},\hat P_t^{\mathrm b}]$ by~\eqref{eq:off_order}, $[\cdot]^+$ is nondecreasing, and $[\cdot]^-$ is nonincreasing, we have $P_t^{\mathrm{ch,d}}\leq\hat P_t^{\mathrm{ch}}\leq P^{\mathrm{ch,max}}$ and $P_t^{\mathrm{dis,d}}\leq\check P_t^{\mathrm{dis}}\leq P^{\mathrm{dis,max}}$. The IT power~\eqref{eq:off_itpower} and the grid power~\eqref{eq:off_grid} are affine in the service rates, the cooling power, and the net BESS power, so the dispatched grid power equals $(1-\beta_t)\check P_t+\beta_t\hat P_t=P_t^{\mathrm d}$.

Next, we verify the time-coupled constraints, which involve the remaining workloads, the hall temperature, and the BESS energy. For the remaining workloads,~\eqref{eq:off_work} gives
\begin{equation}\label{eq:app_work}
 w_{j,t+1}=W_j-\sum_{\tau=a_j}^{t}\mu_{j,\tau}^{\mathrm d}\Delta t
\end{equation}
on the realized trajectory, and likewise for the endpoint trajectories. By~\eqref{eq:off_share}, $\check\mu_{j,\tau}\leq\mu_{j,\tau}^{\mathrm d}\leq\hat\mu_{j,\tau}$, so~\eqref{eq:app_work} implies $\hat w_{j,t}\leq w_{j,t}\leq\check w_{j,t}$ for all $t$. Since $\hat w_{j,d_j+1}=\check w_{j,d_j+1}=0$ by~\eqref{eq:offline_band_complete}, we have $w_{j,d_j+1}=0$, and the realized trajectory satisfies the completion requirement.

For the hall temperature, let $P_t^{\mathrm{heat}}=P_t^{\mathrm{IT}}-\kappa_tP_t^{\mathrm c}$ denote the net heat input and $\nu=1-U\Delta t/C^{\mathrm{th}}$, which is nonnegative because the time step does not exceed the thermal time constant, $\Delta t\leq C^{\mathrm{th}}/U$, for the parameters in Table~\ref{tab:env_parameters}. Then~\eqref{eq:off_thermal} gives
\begin{equation}\label{eq:app_thermal}
 \theta_{t+1}=\nu^{t}\theta_1+\frac{\Delta t}{C^{\mathrm{th}}}\sum_{\tau=1}^{t}\nu^{t-\tau}\big(P_\tau^{\mathrm{heat}}+U\theta_\tau^{\mathrm{amb}}\big)
\end{equation}
on each trajectory. Since $P_t^{\mathrm{heat}}$ is affine in the service rates and the cooling power, its dispatched value is $P_t^{\mathrm{heat,d}}=(1-\beta_t)\check P_t^{\mathrm{heat}}+\beta_t\hat P_t^{\mathrm{heat}}$ and, by~\eqref{eq:off_order}, lies in $[\hat P_t^{\mathrm{heat}},\check P_t^{\mathrm{heat}}]$. Because $\nu^{t-\tau}\geq0$,~\eqref{eq:app_thermal} yields $\hat\theta_t\leq\theta_t\leq\check\theta_t$, and the realized hall temperature satisfies the band in~\eqref{eq:off_tbound}.

For the BESS energy, the charging and discharging powers on each trajectory are $[P_t^{\mathrm b}]^+$ and $[P_t^{\mathrm b}]^-$ under complementarity, so~\eqref{eq:off_battery} gives
\begin{equation}\label{eq:app_soc}
 e_{t+1}=e_1+\sum_{\tau=1}^{t}f^{\mathrm b}\big(P_\tau^{\mathrm b}\big)\Delta t
\end{equation}
with $f^{\mathrm b}(x)=\eta^{\mathrm{ch}}[x]^+-[x]^-/\eta^{\mathrm{dis}}$. The function $f^{\mathrm b}$ is nondecreasing because $\eta^{\mathrm{ch}},\eta^{\mathrm{dis}}>0$, and $P_t^{\mathrm{b,d}}$ lies in $[\check P_t^{\mathrm b},\hat P_t^{\mathrm b}]$, so~\eqref{eq:app_soc} implies $\check e_t\leq e_t\leq\hat e_t$, and the realized BESS energy satisfies~\eqref{eq:off_soc}.

Finally, since $g_t\geq0$ and $P_t^{\mathrm d}\leq\hat P_t$, the realized emissions satisfy $\sum\nolimits_{t\in\mathcal{T}}g_tP_t^{\mathrm d}\Delta t\leq\sum\nolimits_{t\in\mathcal{T}}g_t\hat P_t\Delta t\leq C^{\mathrm{bud}}$ by~\eqref{eq:offline_band_carbon}. Therefore, for any requested power sequence within the intervals, the dispatched controls form an achievable carbon-aware control strategy.
\end{IEEEproof}

\end{document}